\documentclass[12pt,leqno]{article}
\usepackage{amsmath, amsthm, amsfonts, amssymb, color}
\usepackage{mathrsfs}
\theoremstyle{definition}

\newcommand{\scr}[1]{\mathscr #1}
\definecolor{wco}{rgb}{0.5,0.2,0.3}

\numberwithin{equation}{section} \theoremstyle{remark}

\newcommand{\ua}{\uparrow}

\title{{\bf  From Super Poincar\'e   to Weighted Log-Sobolev and Entropy-Cost  Inequalities}
\footnote{Supported in part by NNSFC(10121101) and the 973-Project
in China.} }
\author{
{\bf Feng-Yu Wang}\\
\footnotesize{School of Mathematical Sciences, Beijing Normal
University, Beijing 100875, China}\\
\footnotesize{Present: Department of Mathematics, University of Swansea, Singleton Park, SA2 8PP, UK}}
\begin{document}
\def\R{\mathbb R}  \def\ff{\frac} \def\ss{\sqrt} \def\BB{\mathbb
B}
\def\N{\mathbb N} \def\kk{\kappa} \def\m{{\bf m}}
\def\dd{\delta} \def\DD{\Delta} \def\vv{\varepsilon} \def\rr{\rho}
\def\<{\langle} \def\>{\rangle} \def\GG{\Gamma} \def\gg{\gamma}
  \def\nn{\nabla} \def\pp{\partial} \def\tt{\tilde}
\def\d{\text{\rm{d}}} \def\bb{\beta} \def\aa{\alpha} \def\D{\scr D}
\def\E{\scr E} \def\si{\sigma} \def\ess{\text{\rm{ess}}}
\def\beg{\begin} \def\beq{\begin{equation}}  \def\F{\scr F}
\def\Ric{\text{\rm{Ric}}} \def\Hess{\text{\rm{Hess}}}\def\B{\mathbb B}
\def\e{\text{\rm{e}}} \def\ua{\underline a} \def\OO{\Omega} \def\sE{\scr E}
\def\oo{\omega}     \def\tt{\tilde} \def\Ric{\text{\rm{Ric}}}
\def\cut{\text{\rm{cut}}} \def\P{\mathbb P} \def\ifn{I_n(f^{\bigotimes n})}
\def\C{\scr C}      \def\aaa{\mathbf{r}}     \def\r{r}
\def\gap{\text{\rm{gap}}} \def\prr{\pi_{{\bf m},\varrho}}  \def\r{\mathbf r}
\def\Z{\mathbb Z} \def\vrr{\varrho} \def\ll{\lambda}
\def\L{\scr L}\def\Tt{\tt} \def\TT{\tt} \def\LL{\Lambda}

\maketitle
\begin{abstract} We derive weighted log-Sobolev inequalities from a
class of super Poincar\'e inequalities. As an application, the
Talagrand inequality with larger distances are obtained. In
particular, on a complete connected Riemannian manifold, we prove
that the $\log^\dd$-Sobolev inequality with $\dd\in (1,2)$ implies
the $L^{2/(2-\dd)}$-transportation cost  inequality

$$W^\rr_{2/(2-\dd)}(f\mu,\mu)^{2/(2-\dd)}\le C\mu(f\log f),\ \ \mu(f)=1,
f\ge 0$$ for some constant $C>0$,  and they are equivalent if the
curvature of the corresponding generator is bounded below. Weighted
log-Sobolev and entropy-cost inequalities are also derived for a
large class of probability measures on $\R^d$.
\end{abstract} \noindent
 \vskip 2cm

\section{Introduction}

Let $(E,\rr)$ be a Polish space and $\mu$ a probability measure on
$E$. For $p\ge 1$ we define the $L^p$-Wasserstein distance (or the
$L^p$-transportation cost) by

$$W^\rr_{p}(\mu_1,\mu_2):= \bigg\{\inf_{\pi\in\scr C(\mu_1,\mu_2)}\int_{E\times
E} \rr(x,y)^p \pi(\d x,\d y)\bigg\}^{1/p}$$ for probability measures
$\mu_1,\mu_2$ on $E,$ where $\scr C(\mu_1,\mu_2)$ is the class of
probability measures on $E\times E$ with marginal distributions
$\mu_1$ and $\mu_2.$

According to \cite[Corollary 4]{BV},

$$W^\rr_p(f\mu, \mu)^{2p} \le C\mu(f\log f),\ \ \ f\ge 0,\mu(f)=1$$
holds for some $C>0$ provided $\mu(\e^{\ll
\rr(o,\cdot)^{2p}})<\infty$ for some $\ll>0,$ where $o\in E$ is a
fixed point. See also \cite{DGW} for $p=1.$ Furthermore, it is easy
to derive from \cite[Theorem 1.15]{G06} that for any $q\in [1, 2p)$,
there exists $C>0$ such that

\beq\label{G06}W^\rr_q(f\mu,\mu)^{2p}\le C\mu(f\log f),\ \ \ f\ge
0,\mu(f)=1\end{equation} if and only if  $\mu(\e^{\ll
\rr(o,\cdot)^{2p}})<\infty$ for some $\ll>0.$ In general, however,
this concentration of $\mu$ does not imply (\ref{G06}) for $q=2p.$
Indeed, there exist a plentiful examples where
$\mu(\e^{\ll\rr(o,\cdot)^2})<\infty$ for some $\ll>0$ but there is
no any constant $C>0$ such that the Talagrand  inequality

\beq\label{T2}W^\rr_2(f\mu,\mu)^{2}\le C\mu(f\log f),\ \ \ f\ge
0,\mu(f)=1\end{equation} holds, see e.g. \cite{BLW} for examples
with $\mu(\e^{\ll\rr(o,\cdot)^2})<\infty$ for some $\ll>0$ but the
Poincar\'e inequality does not hold, which is  weaker than
(\ref{T2}) (see \cite[Section 7]{OV} or \cite[Section 4.1]{BGL}).

Therefore, to derive (\ref{G06}) with $q=2p$, one needs something
stronger than the corresponding concentration of $\mu$. In  fact, it
is now well known in the literature that, the Talagrand inequality
follows from the log-Sobolev inequality for a class of local
Dirichlet forms, see \cite{T, OV, BGL, W04, Shao} and references
within.

In this paper, we aim to derive (\ref{G06}) with $q=2p$, i.e.

\beq\label{T2p} W^\rr_{2p}(f\mu,\mu)^{2p}\le C\mu(f\log f),\ \ \ f\ge
0,\mu(f)=1,\end{equation} by using functional inequalities stronger
than the log-Sobolev one.

To this end, in Section 2 we study the weighted log-Sobolev
inequality

$$ \mu(f^2\log f^2)\le C\mu(\aa\circ \rr(o,\cdot) \GG(f,f)),\ \ \
\mu(f^2)=1$$ for a positive function $\aa(r)\to 0$ as $r\to\infty$
and a nice square field $\GG$. Combining this with known results on
log-Sobolev and the Talagrand inequality, we derive (\ref{T2}) with
the original distance $\rr$ replaced by a larger one, which is
induced by the weighted square field $\aa\circ\rr(o,\cdot)\GG$. In
particular, we have the following result on a Riemannian manifold.

Let $M$ be a connected complete Riemannian manifold, and $\mu(\d x)=
\e^{V(x)}\d x$  a probability measure on $M$ for some $V\in C(M).$
We shall use the following super Poincar\'e inequality (see
\cite{W00})

\beq\label{SP} \mu(f^2)\le r\mu(|\nn f|^2) +\bb(r)\mu(|f|)^2,\ \ \
r>0\end{equation} to establish the corresponding weighted
log-Sobolev inequality

\beq\label{WLS}  \mu(f^2\log f^2)\le C\mu(\aa\circ \rr(o,\cdot) |\nn
f|^2),\ \ \ \mu(f^2)=1.\end{equation} By \cite[Theorem 1.1]{W04},
(\ref{WLS}) implies

\beq\label{WT} W_2^{\rr_\aa}(f\mu, \mu)^2\le C\mu(f\log f),\ \ \
f\ge 0, \mu(f^2)=1,\end{equation} where $\rr_\aa$ is the Riemannian
distance induced by the metric

 \beq\label{MM*} \<X,Y\>':= \ff
1{\aa\circ\rr(o,x)}\<X,Y\>,\ \ \ X,Y\in T_x M,\ x\in
M.\end{equation} The main result of the paper is the following.

\beg{thm}\label{T1.1} Assume that $(\ref{SP})$ holds for some
positive decreasing $\bb\in C((0,\infty))$ such that

$$\eta(s):= \big(\log (2s)\big) \big(1\land \bb^{-1} ( s /2)\big),\ \ \ s\ge 1$$ is
bounded, where $\bb^{-1}(s):=\inf\{t\ge 0:\ \bb(t)\le s\}$. Then
 $(\ref{WLS})$ holds for some $C>0$ and

$$\aa(s):= \sup_{t\ge \mu(\rr(o,\cdot)\ge s-2)^{-1}} \eta (t), \ \ \ s\ge 0.$$
 Consequently,  $(\ref{WT})$ holds.
\end{thm}

The following consequences show that the above result is sharp in
specific situations.

\beg{cor} \label{C1.2} Let $\dd\in (1,2).$

$(a)\ (\ref{SP})$ with $\bb(r)= \exp[c(1+r^{-1/\dd}) ]$ implies
$(\ref{WLS})$ with $$\aa(s):= (1
+\rr(o,\cdot))^{-2(\dd-1)/(2-\dd)}$$ and $(\ref{WT})$ with
$\rr_\aa(x,y)$ replaced by
$$\rr(x,y)(1+\rr(o,x)\lor\rr(o,y))^{(\dd-1)/(2-\dd)}.$$
Consequently, it implies

\beq\label{WT'}W^\rr_{2/(2-\dd)}(f\mu,\mu)^{2/(2-\dd)}\le C\mu(f\log
f),\ \ \mu(f)=1, f\ge 0\end{equation}for some constant $C>0.$

$(b)$ If $V\in C^2(M)$ with $\Ric-\Hess_V$ bounded below, then the
following are equivalent to each other: \beg{enumerate} \item[$(1)$]
$(\ref{SP})$  with $\bb(r)= \exp[c(1+r^{-1/\dd}) ]$ for some
constant $c>0;$ \item[$(2)$] $(\ref{WLS})$  with $\aa(s):= (1
+\rr(o,\cdot))^{-2(\dd-1)/(2-\dd)}$  for some $C>0;$ \item[$(3)$]
$(\ref{WT})$ for some $C>0$ and $\rr_\aa(x,y)$ replaced by
$\rr(x,y)(1+\rr(o,x)\lor\rr(o,y))^{(\dd-1)/(2-\dd)};$
\item[$(4)$] $(\ref{WT'})$ for some $C>0;$
\item[$(5)$] $ \mu(\exp[\ll\rr(o,\cdot)^{2/(2-\dd)}])<\infty$
for some $\ll>0.$\end{enumerate} \end{cor}

We remark that (\ref{SP}) with $\bb(r)= \exp[c(1+r^{-1/\dd})]$ for
some $c>0$ is equivalent to the following $\log^\dd$-Sobolev
inequality mentioned in the abstract (see  \cite{W00, W00b, GW,
Wbook} for more general results on (\ref{SP}) and the $F$-Sobolev
inequality)

$$\mu(f^2\log^\dd(1+f^2))\le C_1\mu(|\nn f|^2)+C_2, \ \ \
\mu(f^2)=1.$$  Since due to \cite[Corollary 5.3]{W00b} if (\ref{SP})
holds with $\bb(r)=\exp[c(1+r^{-1/\dd})]$ for some $\dd>2$ then $M$
has to be compact, as a complement to  Corollary \ref{C1.2} we
consider  the critical case  $\dd= 2$ in the next Corollary.

\beg{cor} \label{C1.3} $(\ref{SP})$ with $\bb(r)= \exp[c(1+r^{-1/2})
]$ for some $c>0$ implies $(\ref{WLS})$ with $\aa(s):= \e^{-c_1s}$
for some $c_1>0$ and $(\ref{WT})$ with $\rr_\aa(x,y)$ replaced by

$$\rr(x,y)\e^{c_2[\rr(o,x)\lor\rr(o,y)]}\ge
\e^{c_3\rr(x,y)}-1$$ for some $c_2,c_3>0.$ If $\Ric-\Hess_V$ is
bounded below, they are all equivalent to the concentration
$\mu(\exp[\e^{\ll \rr(o, \cdot)}])<\infty$ for some $\ll>0.$
\end{cor}

\paragraph{Example 1.1.} Let  $\Ric$ be  bounded below.
Let $V\in C(M)$ be such that $V+ a\rr(o,\cdot)^\theta$ is bounded
for some $a>0$ and $\theta\ge 2.$ By \cite[Corollaries 2.5 and
3.3]{W00}, (\ref{SP}) holds for $\dd= 2(\theta-1)/\theta$. Then
Corollary \ref{C1.2} implies

$$ W^\rr_{\theta}(f\mu,\mu)^\theta \le C\mu(f\log f), \ \ \ f\ge 0, \mu(f)=1$$ for
some constant $C>0.$

In this inequality $\theta$ could not be replaced by any larger
number, since $W^\rr_\theta\ge W^\rr_1$ and by Proposition \ref{P2.3} below
for any $p\ge 1$ the inequality

$$W^\rr_1(f\mu,\mu)^p \le C\mu(f\log f), \ \ \ f\ge 0, \mu(f)=1$$
implies $\mu(\e^{\ll \rr(o,\cdot)^p})<\infty$ for some $\ll>0,$
which fails when $p>\theta$ for $\mu$ specified above.

\
\paragraph{Example 1.2.} In the situation of Example 1.1 but let $V+
\exp[\si\rr(o,\cdot)]$ be bounded for some $\si>0.$ Then by
\cite[Corollaries 2.5 and 3.3]{W00}, (\ref{SP}) holds with $\bb(r)=
\exp[c(1+r^{-1/2}) ]$ for some $c>0$. Hence, by Corollary
\ref{C1.3},

\beq\label{N1.9}\inf_{\pi\in\scr C(\mu,f\mu)}\int_{M\times M}
\rr(x,y)^2\e^{c_1\rr(x,y)} \pi(\d x,\d y)\le C\mu(f\log f),\ \ \
f\ge 0, \mu(f)=1\end{equation} holds for some $c_1, C>0.$

On the other hand, it is easy to see from Jensen's inequality that
the left hand side is larger than

$$(\exp[c_2W^\rr_1(\mu,f\mu)]-1)^2$$ for some $c_2>0.$
So, by Proposition \ref{P2.3} below  (\ref{N1.9}) implies
$\mu(\exp[\exp(\ll \rr(o,\cdot))])<\infty$ holds for any $\ll>0,$
which is the exact concentration property of the given measure
$\mu$.

\

 In the next section we study the super Poincar\'e and the weighted
log-Sobolev inequality in an abstract framework, and complete proofs
of the above results are presented in Section 3. Finally, weighted
log-Sobolev and transportation cost inequalities are also studied
for  probability measures on $\R^d$ by using concentrations.

\section{From super Poincar\'e to weighted log-Sobolev inequalities}

We shall work with a diffusion framework as in \cite{BLW}. Let
$(E,\F,\mu)$ be a separable complete probability space, and let
$(\E,\D(\E))$ be a conservative symmetric local Dirichlet form on
$L^2(\mu)$ with domain $\D(\E)$ in the following sense.  Let $\scr
A$ be a dense subspace of $\D(\E)$ under the $\E_1^{1/2}$-norm ($
\E_1(f,f) = \| f\| _2^2 + \E (f,f)$) 
 which is composed of bounded functions, stable under products  and composition with Lipschitz
functions on $\R$. Let $\GG: \scr A\times \scr A\to \scr M_b$ be a
bilinear mapping, where $\scr M_b$ is the set of all bounded
measurable functions on $E$, such that

\

(1) $\GG(f,f)\ge 0$ and $\E(f,g)=\mu(\GG(f,g))$ for $f,g\in \scr A$;

(2) $\GG(\phi\circ f, g)= \phi'(f) \GG(f,g)$ for $f,g\in \scr A$ and
$\phi\in C_b^\infty(\R);$

(3) $\GG(fg,h)= g\GG(f,h)+ f\GG(g,h)$ for $f,g,h\in \scr A$ with
$fg\in \scr A.$

\

\noindent It is easy to see that the positivity and the bilinear
property imply $\GG(f,g)^2\le \GG(f,f)\GG(g,g)$ for all  $f,g\in
\scr A$. For simplicity we set below $\GG (f,f) = \GG(f)$ and
$\E(f,f) = \E(f)$.

We shall denote by $\scr A_{\rm { loc}}$ the set of functions $f$
such that for any integer $n$, the truncated function $f_{n    }=
\min(n, \max(f,-n))$ is in $\scr A$. For such functions, the
bilinear map $\GG $ automatically extends and shares the same
properties than for functions in $\scr A$. 

Next, let \def\vrr{\varrho} $\vrr\in \scr A_{\rm{loc}}$ be positive
such that $\GG(\vrr,\vrr)\le 1.$ We shall start from the super
Poincar\'e inequality

\beq\label{SP'} \mu(f^2)\le r\E(f,f)+\bb(r)\mu(|f|)^2,\ \ \
r>0.\end{equation} To derive the desired weighted log-Sobolev
inequality

\beq\label {WLS'} \mu(f^2\log f^2)\le C\mu(\GG(f,f) \aa\circ\vrr),\
\ \ \mu(f^2)=1,\end{equation} we shall also need the following
Poincar\'e inequality

\beq\label{P} \mu(f^2)\le C_0\E(f,f) +\mu(f)^2\end{equation} for
some $C_0>0$. Here and in what follows, the reference function $f$
is taken from $\scr A.$

\beg{thm} \label{T2.1} Assume $(\ref{P})$ holds for some $C_0>0.$
Then $(\ref{SP'})$ implies  $(\ref{WLS'})$ for some constant $C>0$
and $\aa$ given in Theorem $\ref{T1.1}.$
\end{thm}

\beg{proof} (a) Let $\Phi(s)=\mu(\vrr\ge s)$ which decreases to zero
as $s\to\infty.$ We may take $r_0>0$ such that

\beq\label{2.1} r_0(1+ \sup_{s\ge 1}\eta(s))\le \ff 1
{32}\end{equation} and

\beq\label{2.2} \bb^{-1} (\e^{r_0^{-1}}/4)\le 1.\end{equation} For a
fixed number $r\in (0,r_0]$ we define

\beg{equation*}\beg{split} &h_n= \big((\vrr
-\Phi^{-1}(2\e^{-r^{-1}})-n)_+\land 1\big)\big((n+2
+\Phi^{-1}(2\e^{-r^{-1}})-\vrr)_+\land 1\big),\\
&\dd_n= \Big(\log \ff 2 {\Phi(n+ \Phi^{-1} (2\e^{-r^{-1}}))}\Big)
\bb^{-1} \Big(\ff 1 {2 \Phi(n+ \Phi^{-1} (2\e^{-r^{-1}}))}\Big),\\
& B_n= \{n\le \vrr -\Phi^{-1} (2\e^{-r^{-1}})\le n+2\},\ \ \ \ n\ge
0.\end{split}\end{equation*} Then

\beq\label{2} \sum_{n=0}^\infty h_n^2\ge \ff 1 2 1_{\{\rr\ge 1+
\Phi^{-1} (2\e^{-r^{-1}})\}}.\end{equation} By (\ref{SP'}) and
noting that

$$\mu(|f|h_n)^2\le \mu(f^2h_n^2) \mu(\vrr> n+\Phi^{-1}
(2\e^{-r^{-1}}))\le\mu(f^2h_n^2) \Phi(n+\Phi^{-1}
(2\e^{-r^{-1}})),$$ we have

\beg{equation*} \beg{split}&\sum_{n=0}^\infty \mu(f^2 h_n^2) \le
\sum_{n=0}^\infty \Big\{  r_n \mu\big(\GG(fh_n,fh_n)\big)
+\bb(r_n)\mu(|f|h_n)^2\Big\}\\
&\le \sum_{n=0}^\infty \Big\{ \ff{2 r_n}{\dd_n}\mu(\GG(f,f) \dd_n
1_{B_n}) + 2 r_n \mu(f^2 1_{B_n})+\bb(r_n) \Phi(n+\Phi^{-1}
(2\e^{-r^{-1}})) \mu(f^2h_n^2)\Big\}\end{split}\end{equation*} for
$r_n>0.$ Since by (\ref{2.2}) and the definition of $\aa$

$$\aa(s) \ge \dd_n \  \text{for}\ s \ge n+ 2 +\Phi^{-1} (2\e^{-r^{-1}}),$$
letting $r_n= \dd_n r$ we obtain

\beq\label{3} \beg{split}\sum_{n=0}^\infty \mu(f^2 h_n^2) \le
\sum_{n=0}^\infty \Big\{ & 2r\mu(\GG(f,f) \aa\circ\vrr
1_{B_n}) + 2r \dd_n \mu(f^2 1_{B_n})\\
& +\bb(r\dd_n) \Phi(n+\Phi^{-1} (2\e^{-r^{-1}}))
\mu(f^2h_n^2)\Big\}.\end{split}\end{equation} Noting that

$$A:=r\log \ff 2 {\Phi(n+\Phi^{-1} (2\e^{-r^{-1}}))} \ge r\log \ff 2 {\Phi(\Phi^{-1}
(2\e^{-r^{-1}}))}=1,$$ we have

$$\bb(\dd_n r) = \bb\Big(A \bb^{-1}\Big( \ff 1 {2 \Phi(n+\Phi^{-1}
(2\e^{-r^{-1}}))}\Big)\Big)\le \ff 1{2\Phi(n+\Phi^{-1}
(2\e^{-r^{-1}}))}.$$ Thus, by (\ref{3}) and (\ref{2.1}) and the fact
that $\dd_n\le \sup\eta,$  we arrive at

$$\sum_{n=0}^\infty \mu(f^2 h_n^2) \le \sum_{n=0}^\infty \Big\{
2r\mu(\GG(f,f) \aa\circ\vrr 1_{B_n}) + \ff 1 8 \mu(f^2) +\ff 1
2\sum_{n=0}^\infty \mu(f^2 h_n^2).$$ It follows from this and
(\ref{2}) that

\beq\label{4} \mu(f^21_{\{\vrr \ge 1 +\Phi^{-1}(2\e^{-r^{-1}})\}})
\le 8 r \mu(\GG(f,f) \aa\circ\vrr) + \ff 1 2 \mu(f^2).\end{equation}

(b) On the other hand, since $\aa$ is decreasing

\beg{equation*}\beg{split} &\mu(f^21_{\{\vrr\le
1+\Phi^{-1}(2\e^{-r^{-1}})\}})\le \mu(f^2
\{(2+\Phi^{-1}(2\e^{-r^{-1}})-\vrr)_+^2\land 1\})\\
&\le 2s \mu(\GG(f,f)1_{\{\vrr\le2+\Phi^{-1}(2\e^{-r^{-1}})\}})+ 2s
\mu(f^2) + \bb(s) \mu(|f|)^2\\
&\le
\ff{2s}{\aa(2+\Phi^{-1}(2\e^{-r^{-1}}))}\mu(\GG(f,f)\aa\circ\vrr)+
2s \mu(f^2) + \bb(s) \mu(|f|)^2,\ \  \
s>0.\end{split}\end{equation*} Taking

$$s= r\aa(2+\Phi^{-1}(2\e^{-r^{-1}}))\le \ff 1 {32}$$
due to (\ref{2.1}), we obtain

$$\mu(f^21_{\{\vrr\le
1+\Phi^{-1}(2\e^{-r^{-1}})\}}) \le 2 r \mu(\GG(f,f) \aa\circ\vrr)
+\ff 1 {16} \mu(f^2) +\bb
\big(r\aa(2+\Phi^{-1}(2\e^{-r^{-1}}))\big)\mu(|f|)^2.$$ Since by
(\ref{2.2}) and the definition of $\aa$

\beg{equation*}\beg{split}
r\aa\big(2+\Phi^{-1}(2\e^{-r^{-1}})\big)&\ge  \Big(r\log
\ff{2}{\Phi(\Phi^{-1}(2\e^{-r^{-1}}))}\Big)
\bb^{-1}\Big(\ff 1 {2\Phi(\Phi^{-1}(2\e^{-r^{-1}}))}\Big)\\
&=  \bb^{-1}\Big(\ff
{\e^{r^{-1}}}{4}\Big),\end{split}\end{equation*} we obtain

$$\mu(f^21_{\{\vrr\le 1+\Phi^{-1}(2\e^{-r^{-1}})\}}) \le 2 r
\mu(\GG(f,f) \aa\circ\vrr) +\ff 1 {16} \mu(f^2) +\ff {\e^{r^{-1}}}4
\mu(|f|)^2.$$ Combining this with (\ref{4}) we conclude that

$$\mu(f^2)\le 40 r\mu(\GG(f,f) \aa\circ\vrr) +\e^{r^{-1}}\mu(|f|)^2,\ \ \ r\in (0,r_0].$$
Therefore, there exists a constant $c>0$ such that

\beq\label{**}\mu(f^2)\le r\mu(\GG(f,f) \aa\circ\vrr) +
\e^{c(1+r^{-1})}\mu(|f|)^2,\ \ \ r>0.\end{equation} According to
e.g. \cite[Corollary 1.3]{W00b}, this is equivalent to the defective
weighted log-Sobolev inequality

\beq\label{5}\mu(f^2\log f^2)\le C_1\mu(\GG(f,f) \aa\circ\vrr)+
C_2,\ \ \ \mu(f^2)=1.\end{equation}

(c) Finally, for any $f$ with $\mu(f)=0$, it follows from (\ref{P})
that

\beg{equation*}\beg{split} \mu(f^2)&\le \mu(f^2\{(1+R-\vrr)_+^2\land
1\}) +\|f\|_\infty^2 \mu(\vrr\ge R)\\
&\le 2C_0 \mu(\GG (f,f)1_{\{\vrr\le 1+R\}}) + (2C_0+1)\|f\|_\infty^2
\mu(\vrr\ge R)+\mu(f\{(\vrr -R)_+\land 1\})^2\\
&\le  \ff{2C_0}{\aa(1+R)} \mu(\GG (f,f)\aa\circ\vrr) +
2(C_0+1)\|f\|_\infty^2 \mu(\vrr\ge R), \ \ \
R>0.\end{split}\end{equation*} Since $\mu(\vrr\ge R)\to 0$ as
$R\to\infty$, the weighted weak Poincar\'e inequality

$$\mu(f^2)\le \tt\bb(r)\mu(\GG(f,f)\aa\circ\vrr)+r \|f\|_\infty^2,\
\ r>0, \mu(f)=0$$ holds for some positive function $\tt\bb$ on
$(0,\infty).$ By \cite[Propsosition 1.3]{RW01}, this and (\ref{**})
implies the weighted Poincar\'e inequality

$$\mu(f^2)\le C'\mu(\GG(f,f)\aa\circ\vrr)+\mu(f)^2$$ for some
constant $C'>0.$ Combining this with (\ref{5}) we obtain the desired
weighted log-Sobolev inequality (\ref{WLS'}). \end{proof}

\section{Proofs of Theorem \ref{T1.1} and Corollaries}

\emph{Proof of Theorem \ref{T1.1}.} Since $\aa$ is bounded, the
completeness of the original metric implies that of the weighted one
given by (\ref{MM*}).  So, (\ref{WT}) follows from (\ref{WLS}) due to
\cite[Theorem 1.1]{W04} with $p\to 2$. Thus, by Theorem \ref{T2.1}
with $E=M$ and  $\GG(f,f)=|\nn f|^2$, it suffices to prove that
(\ref{SP}) implies the Poincar\'e inequality (\ref{P}) for some
$C_0>0.$ Due to \cite{W00} the super Poincar\'e inequality
(\ref{SP}) implies that the spectrum of $L$ is discrete. Moreover,
since $M$ is connected, the corresponding Dirichlet form is
irreducible so that $0$ is a simple eigenvalue.   Therefore, $L$
possesses a  spectral gap, which is equivalent to the desired
Poincar\'e inequality. \qed

\ \newline To complete the proof of Corollary \ref{C1.2}, in the
spirit of \cite{Ma, BG} we introduce below a deviation inequality
induced by the $L^1$-transportation cost inequality.

\beg{prp} \label{P2.3} Let $\tt\rr: M\times M\to [0,\infty)$ be
measurable. For any $r>0$ and measurable set $A\subset M$ with
$\mu(A)>0$, let

$$A_r=\{x\in M: \tt\rr(x,y)\ge r\ \text{for\ some\ }y\in A\},\ \ \ r>0.$$ If

\beq\label{3.2} W^{\tt\rr}_1(f\mu,\mu)\le \Phi\circ\mu(f\log f),\ \
\ f\ge 0, \mu(f)=1\end{equation} holds for some positive increasing
$\Phi\in C([0,\infty))$, then

\beq\label{3.3} \mu(A_r) \le \exp\big[-\Phi^{-1}(r-\Phi\circ\log
\mu(A)^{-1})\big],\ \ \ r>\Phi\circ\log \mu(A)^{-1},\end{equation}
where $\Phi^{-1}(r):= \inf\{s\ge 0:\ \Phi(s)\ge r\},\ r\ge
0.$\end{prp}

\beg{proof} It suffices to prove for $\mu(A_r)>0.$ In this case,
letting $\mu_A=\mu(\cdot\cap A)/\mu(A)$ and $\mu_{A_r}=\mu(\cdot\cap
A_r)/\mu(A_r)$, we obtain from (\ref{3.2}) that

$$r\le W^{\tt\rr}_1(\mu_A,\mu_{A_r})\le
W^{\tt\rr}_1(\mu_A,\mu)+W^{\tt\rr}_1(\mu_{A_r},\mu)\le \Phi\circ\log
\mu(A)^{-1} +\Phi\circ\log \mu(A_r)^{-1}.$$ This completes the
proof.\end{proof}

\ \newline\emph{Proof of Corollary \ref{C1.2}.} (a) Let $\bb(r) =
\e^{c(1+r^{-1/\dd})}$ for some $c>0$ and $\dd>1.$ It is easy to see
that

$$1\land \bb^{-1}(s/2)\le c_1 \log^{-\dd}(2s),\ \ \ s\ge 1$$ holds for some constant $c_1>0.$ Next,
by \cite[Corollary 5.3]{W00b}, (\ref{SP}) with this specific
function $\bb$ implies

$$\mu(\rr(o,\cdot)\ge s-2)\le c_2\exp[-c_3 s^{2/(2-\dd)}],\ \ \ s\ge 0$$
for some constants $c_2, c_3>0.$ Therefore,

\beq\label{N1} \aa(s)\le c_4(1+s)^{-2(\dd-1)/(2-\dd)},\ \ \ s\ge
0\end{equation}  holds for some constant $c_4>0.$

On the other hand,  for any $x_1,x_2\in M$ let $i\in \{1,2\}$ such
that $\rr(o,x_i)=\rr(o,x_1)\lor \rr(o,x_2).$ Define

$$f(x)= \big(\rr(x,x_i)\land \ff{\rr(o,x_i)} 2 \big)(1+\rr(o,x_i))^{(\dd-1)/(2-\dd)},\ \ \ x\in \R^d.$$
Then

\beg{equation*}\beg{split}&\aa\circ\rr(o,\cdot)|\nn f|^2\le
c_4(1+\rr(o,\cdot))^{-2(\dd-1)/(2-\dd)}|\nn f|^2\\
&\le c_4 1_{\{\rr(o,x_i)/2\le \rr(o,\cdot)\le 3\rr(o,x_i)/2\}}
(1+\rr(o,\cdot))^{-2(\dd-1)/(2-\dd)}
(1+\rr(o,x_i))^{2(\dd-1)/(2-\dd)} \le c_5\end{split}\end{equation*}
for some constant $c_5>0.$ Since by the triangle inequality
$\rr(o,x_i)\ge \ff 1 2 \rr(x_1,x_2)$, this implies that the
intrinsic distance $\rr_\aa$ satisfies

\beg{equation*}\beg{split} &\rr_\aa (x_1,x_2)^2 \ge \ff
{|f(x_1)-f(x_2)|^2} {c_5} \\
&\ge c_6 \rr(x_1,x_2)^2(1+\rr(o,x_1)\lor
\rr(o,x_2))^{2(\dd-1)/(2-\dd)}\ge
c_7\rr(x_1,x_2)^{2/(2-\dd)}\end{split}\end{equation*} for some
constant $c_6,c_7>0.$ Hence the proof of (a) is completed by Theorem
\ref{T1.1}.

(b) Now, assume that

$$\Ric-\Hess_V\ge -K$$ for some $K\ge 0.$ By (a) and Proposition \ref{P2.3}, which ensures the implication
from (4) to (5), it suffices to deduce (1) from (5).  Let

$$h(r)= \mu(\e^{r\rr(o,\cdot)^2}),\ \ \ \ r>0.$$ By \cite[Theorem
5.7]{W00b}, the super Poincar\'e inequality (\ref{SP}) holds with
\beq\label{BB} \bb(r):= c_0\inf_{0< r_1<r} r_1 \inf_{s>0} \ff 1 s
h(2K+12 s^{-1})\e^{s/r_1 -1},\ \ \ r>0\end{equation} for some
constant $c_0>0.$ Since for any $\ll>0$ there exists $c(\ll)>0$ such
that

$$rt^2\le \ll t^{2/(2-\dd)} + c(\ll) r^{1/(\dd-1)},\ \ \ r>0,$$ it
follows from (5)  that

$$h(r)\le c_1 \exp[c_1 r^{1/(\dd-1)}],\ \ \ r>0$$
for some constants $c_1>0.$  Therefore,

$$\bb(r) \le c_2\inf_{0< r_1<r} r_1 \inf_{s>0} \ff 1 s
\exp[c_2s^{-1/(\dd-1)}+ s/r_1],\ \ \ r>0$$ for some $c_2>0.$ Taking
$s= r^{(\dd-1)/\dd}$ and $r_1=r,$ we conclude that

$$ \bb(r)\le \e^{c(1+ r^{-1/\dd})},\ \ \ r>0$$
for some $c>0.$  Thus, (1) holds. \qed

\ \newline\emph{Proof of Corollary \ref{C1.3}.} The proof  is
similar to that of Corollary \ref{C1.2} by noting that (\ref{SP})
with $\bb(r)=\exp[c(1+r^{-1/2})]$ implies $\mu(\rr(o,\cdot)\ge s)\le
\exp[-c\e^{c_1 s}]$ for some $c_1>0,$ see \cite[Corollary
5.3]{W00b}.\qed

\section{Weighted log-Sobolev and transportation cost inequalities  on $\R^d$}

Our main purpose of this section  is to establish the weighted
log-Sobolev inequality for an arbitrary probability measure using
the concentration of this measure. We shall also prove the
 HWI inequality
introduced in \cite{BGL} for the corresponding weighted Dirichlet form.
The main point  is to find square fields (resp. cost
functions) for a given probability measure to satisfy the
log-Sobolev inequality (resp. the Talagrand transportation cost
inequality). So, the line of our study is exactly opposed to existed
references in the literature, see e.g. \cite{GGM, GGM2, Go} and
references within, which provided conditions on the reference
measure such that the log-Sobolev (resp. transportation cost)
inequality holds for a given square field (resp. the corresponding
cost function).

The basic idea of the study comes from Caffarelli \cite{Ca} which
says that for any probability measure  $\mu(\d x):=\e^{V(x)}\d x$
 on $\R^d$, there exists a convex function $\psi$ on $\R^d$ such
 that $\nn\psi$  pushes $\mu$ forward to the standard Gaussian
 measure $\gamma$; that is, letting

 $$y(x):=\nn \psi(x),\ \ \ x\in \R^d,$$
which is one-to-one, one has $\gamma= \mu\circ y^{-1}.$ Furthermore,
$\nn \psi$ is uniquely determined and $\Hess_\psi$ is non-degenerate
with

$$\text{det} (\Hess_\psi) =(2\pi)^{d/2}\e^{V+|\nn \psi|^2/2}.$$
Let

$$\rr(x_1,x_2):= |y(x_1)-y(x_2)|,\ \ \ x_1,x_2\in \R^d.$$
Let $W_2$ be the $L^2$-Wasserstein distance induced by the usual
Euclidiean metric. Due to Talagrand \cite{T}

\beq\label{4.1.1} W_2(\gamma, f^2\gamma)^2\le 2 \gamma(f^2\log
f^2),\ \ \ \ \gamma(f^2)=1.\end{equation}
 Since $\pi\in \scr C
(\mu\circ y^{-1}, (f^2\circ y^{-1})\mu\circ y^{-1})$ if and only if
$\pi\circ (y\otimes y)\in \scr C (\mu, f^2\mu)$, we obtain from
(\ref{4.1.1})  and   the change of variables theorem that

$$W_2^\rr(\mu, f^2\mu)^2 = W_2 (\gamma, (f^2\circ y^{-1})\gamma)^2\le
2\gamma(f^2\circ y^{-1}\log f^2\circ y^{-1})=2\mu(f^2\log f^2),\ \ \
\mu(f^2)=1.$$ Similarly, since

$$\nn (f\circ y^{-1}) = (D y^{-1})(\nn f)\circ y^{-1}= [(D y)\circ
y^{-1}]^{-1} (\nn f)\circ y^{-1}= [(\Hess_\psi)^{-1}\nn f]\circ
y^{-1},$$ where $D y:= (\pp_i y_j)_{d\times d}$, by Gross'
log-Sobolev inequality for $\gamma$ (see \cite{Gross}) we obtain

$$\mu(f^2\log f^2) \le 2 \mu(|(\Hess_\psi)^{-1}\nn f|^2),\ \ \
f\in C_0^\infty(\R^d), \mu(f^2)=1.$$

On the other hand, however, since the transportation $\nn \psi$ is
normally inexplicit, it is hard to estimate the distance $\rr$ and
the matrix $\Hess_\psi$. So, to derive transportation and
log-Sobolev inequalities with explicit distances and Dirichlet
forms, we shall construct, instead of $\nn \psi$, an explicit map
using the concentration of $\mu$, which transports the measure
into the standard Gaussian measure with a perturbation. In many
cases this perturbation is bounded and hence, does not make much
trouble to derive the desired inequalities.
\def\S{\mathbb S}
\subsection{Main results}

In this subsection we provide an explicit  positive function $\aa$
and an explicit distance $\rr$ on  $\R^d$ such that the
log-Sobolev inequality

\beq\label{4.1.2} \mu(f^2\log f^2)\le 2 \mu(\aa|\nn f|^2),\ \ \ f\in
C_0^\infty(\R^d),\ \mu(f^2)=1\end{equation} and the
transportation-cost inequality

\beq\label{4.1.3} W^\rr_2 (\mu, f^2\mu)^2\le 2 \mu(f^2\log f^2),\ \
\ \ \mu(f^2)=1\end{equation} hold. In a special case, we are also
able to present the HWI inequality stronger than (\ref{4.1.2}).

Let us first consider a probability measure $\mu(\d x):=
\e^{V(x)}\d x$ on $[\dd,\infty)$ for some  $\dd\in
[-\infty,\infty)$, where $[-\infty,\infty)$ is regarded as $\R$.
Let

$$\Phi_\dd (r):= \ff 1 {c_\dd}\int_{\dd}^r \e^{-s^2/2}\d s,\ \ \varphi(r):=
\mu([\dd,r))=\int_{\dd}^r \e^{V(x)}\d x,\ \ r\ge \dd,$$ where
$c_\dd:=\int_\dd^\infty \e^{-x^2/2}\d x$ is the normalization.

\beg{thm}\label{T4.1.1} Let $\mu(\d x):=1_{[\dd,\infty)}(x)\e^{V(x)}\d
x$ be a probability measure on $[\dd,\infty).$ For the above defined
$\Phi_\dd$ and $\varphi$, $(\ref{4.1.2})$ and $(\ref{4.1.3})$ hold with
$\R^d$ replaced by $[\dd,\infty)$  for

\beg{equation*}\beg{split} &\aa:= \Big(\ff{\Phi_\dd'\circ
\Phi_\dd^{-1}\circ\varphi}{\varphi'}\Big)^2,\\
&\rr(x,y):= |\Phi_\dd^{-1}\circ\varphi(x)
-\Phi_\dd^{-1}\circ\varphi(y)|,\ \ \ x,y\ge
\dd.\end{split}\end{equation*}  Furthermore,

\beq\label{4.1.4} \mu(f^2\log f^2) +W_2^\rr(\mu,f^2\mu)^2 \le
2\ss{2\mu(\aa {f'}^2)} W_2^\rr(\mu, f^2\mu),\ \ \ f\in
C_0^\infty([\dd,\infty)), \mu(f^2)=1.\end{equation} \end{thm}

The inequality (\ref{4.1.4}), linking the Wasserstein distance, the
relative entropy  and the energy, is called the HWI inequality in
\cite{BGL} and \cite{OV2}.

To extend this result to $\R^d$ for $d\ge 2$, we consider the
polar coordinate $(r,\theta)\in [0,\infty)\times \S^{d-1}$, where
$\S^{d-1}$ is the unit sphere in $\R^d$ with the induced metric.
Then $\mu$ can be represented as

$$\d\mu= c(d) r^{d-1} \e^{V(r\theta)}\d r\d \theta =:G(r,\theta)\d r\d\theta,$$
where $\d\theta$ is the normalized volume measure on $\S^{d-1}$, and
$c(d)/d$ equals to the volume of the unit ball in $\R^d$. Let
$B(0,r):=\{x\in\R^d: |x|< r\}$ and

 \beg{equation*}\beg{split}
&\Phi_0(r):= \int_{B(0,r)}
\ff{\e^{-|x|^2/2}\d x}{(2\pi)^{d/2}},\ \ r\ge 0,\\
&h(\theta):=\int_0^\infty s^{d-1}\e^{V(s\theta)}\d s,\ \ \
\theta\in
\S^{d-1},\\
&\varphi_\theta(r):= \ff 1 {h(\theta)}\int_0^r
s^{d-1}\e^{V(s\theta)}\d s,\ \ \ \theta\in\S^{d-1}, r\ge
0.\end{split}\end{equation*} Since $\mu(\R^d)=1$, we have
$h(\theta)\in (0,\infty)$ for a.e. $\theta\in \S^{d-1}$.

We shall prove that the map

$$x\mapsto \Phi_0^{-1}\circ\varphi_{\ff{x}{|x|}}(|x|)\ff x {|x|}$$
transports $\mu$ into a Gaussian measure with density
$h\circ\theta$. Thus, to derive the desired inequalities for
$\mu$, we  need a regularity property of this transportation
specified in the following result.

\beg{thm} \label{T4.1.2} Let  $r(x):=|x|, \theta(x):= \ff {x}{|x|},
x\in \R^d$. If $C(h):=\sup_{\theta_1,\theta_2\in \S^{d-1}}\ff
{h(\theta_1)}{h(\theta_2)}<\infty,$ then $(\ref{4.1.3})$ holds for

$$\rr(x_1,x_2):= C(h)^{-1/2}|(\Phi_0^{-1}\circ\varphi_{\theta}(r)\theta)(x_1)
- (\Phi_0^{-1}\circ\varphi_{\theta}(r)\theta)(x_2)|,\ \ \
x_1,x_2\in \R^d.$$ If moreover $\varphi_\theta(r)$ is
differentiable in $\theta$ then $(\ref{4.1.2})$ holds for

$$\aa:=C(h) \inf_{\vv>0}\max\Big\{\ff
{(1+\vv)r^2}{(\Phi_0^{-1}\circ\varphi_\theta (r))^2},\
\ff{(\Phi_0'\circ\Phi_0^{-1}\circ\varphi_\theta(r))^2}{({\varphi_\theta}'(r))^2}
+\ff{(1+\vv^{-1})|\nn_\theta \varphi_\theta
(r)|^2}{({\varphi_\theta}'(r)\Phi_0^{-1}\circ\varphi_\theta(r))^2}\Big\}.
$$  If, in
particular, $h$ is constant (it is the case if $V(x)$ depends only
on $|x|$), then the following HWI inequality holds:

\beq\label{4.1.7} \mu(f^2\log f^2) + W_2^\rr(\mu, f^2\mu)^2 \le
2\ss{2\mu(\aa|\nn f|^2)}\, W_2^\rr(\mu,f^2\mu),\ \ \ f\in
C_0^\infty(\R^d), \mu(f^2)=1,\end{equation} for

$$\aa:=\max\Big\{\ff
{r^2}{(\Phi_0^{-1}\circ\varphi (r))^2},\
\ff{(\Phi_0'\circ\Phi_0^{-1}\circ\varphi(r))^2}{(\varphi'(r))^2}\Big\}$$
and $\varphi=\varphi_\theta$ is independent of $\theta.$
\end{thm}

Note that if $V$ is locally bounded and $\zeta(r):= \sup_{|x|=r} V(x)$
satisfies $\int_0^\infty r^{d-1}\e^{\zeta(r)}\d r<\infty,$ then
$C(h)<\infty.$ Thus, Theorem \ref{T4.1.2} applies to a large number
of probability measures. In particular, we have the following
concrete result.

\beg{cor}\label{C4.1.3} Let $V$ be differentiable such that $\mu(\d
x):=\e^{V(x)}\ d x$ is a probability measure and

\beq\label{V} -c_1 |x|^{\dd-1} \le \<\nn V(x), \nn |x|\> \le - c_2
|x|^{\dd-1}\end{equation} holds for some constants $\dd, c_1, c_2>0$
and large $|x|$. If there exists a constant $c_3>0$ such that

\beq\label{V2} |\nn_\theta V|\le c_3,\end{equation} where
$\nn_\theta$ is the gradient on $\S^{d-1}$ at point $\theta,$ then
there exists a constant $c>0$ such that

\beq\label{4.1.9} \mu(f^2\log f^2)\le c\mu((1+|\cdot|)^{2-\dd}|\nn
f|^2),\ \ \ f\in C_0^\infty(\R^d), \mu(f^2)=1.\end{equation}
Consequently,

\beq\label{4.1.10} W_2^{\tt\rr}(\mu,f^2\mu)^2\le c'\mu(f^2\log f^2),\
\ \ \ \mu(f^2)=1\end{equation} holds for some constant $c'>0$ and
$$\tt\rr(x,y):=
\ff{|x-y|}{(1+|x|\lor |y|)^{1-\dd/2}},\ \ \ x,y\in \R^d.$$
\end{cor}

\paragraph{Remark.} (a) The inequalities presented in Corollary
\ref{C4.1.3} are sharp in the sense that (\ref{4.1.10}) (and hence
also (\ref{4.1.9})) implies $\mu(\e^{\ll r^\dd})<\infty$ for some
$\ll>0$, which is the exact concentration of $\mu$. This follows
from \cite[Corollary 3.2]{BG} and the fact that
$\tt\rr(0,x)\thickapprox |x|^{\dd/2}$ for large $|x|$.

(b) When $V$ is strictly concave, the matrix

$$\LL[v_1,v_2]:= \int_0^1 s (-\Hess_V)((1-s)v_1+s  v_2)\d s$$
is strictly positive definite for any $v_1,v_2\in\R^d$. It is
proved by Kolesnikov (see \cite[Corollary 3.1]{K}) that

\beq\label{4.1.11} \mu(f^2\log f^2) \le \int_{\R^d} \<\LL[T_f,
\cdot]^{-1} \nn f,\nn f\>\d\mu,\ \ \ f\in C_0^\infty(\R^d),
\mu(f^2)=1,\end{equation} where $x\mapsto T_f(x)$ is the optimal
transport of $f^2\mu$ to $\mu$. In particular, for $V(x):=
-|x|^\dd+c$ with $\dd >2$ and a constant $c$, \cite[Example
3.2]{K} implies (\ref{4.1.9}) for even smooth function $f^2$. But
Corollary \ref{C4.1.3} works for more general $V$ and  all smooth
function $f$.

(c) Recently, Gentil, Guillin and Miclo \cite{GGM} (see \cite{GGM2,
Go} for further study) established a Talagrand type inequality for
$V(x)= -|x|^\dd +c$ with $\dd\in [1,2]$ and a constant $c$.
Precisely, there exist constants $a, D>0$ such that

\beq\label{4.1.12}\inf_{\pi\in \scr C(\mu,f^2\mu)} \int_{\R^d\times
\R^d} L_{a,D}(x-y)\pi(\d x,\d y) \le D \mu(f^2\log f^2),\ \ \ \
\mu(f^2)=1,\end{equation} where

$$L_{a,D}(x):= \beg{cases} \ff{|x|^2} 2, &\text{if}\ |x|\le a,\\
\ff{a^{2-\dd}}\dd |x|^\dd + \ff{a^2(\dd-2)}{2\dd},
&\text{otherwise.}\end{cases}$$ Since $L_{a,D}(x-y)\ge \vv
\tt\rr(x,y)^2$ for some constant $\vv>0$,  this inequality implies
(\ref{4.1.10}) for $\dd\in [1,2].$ But (\ref{4.1.12}) is yet unavailable
for $\dd\notin [1,2]$ while (\ref{4.1.10}) holds for more general $V$.
In particular, if $\dd>2$ then (\ref{4.1.10}) with $\tt\rr(x,y)\ge
c(|x-y|\lor |x-y|^{\dd/2})$ for some $c>0$, which is new in the
literature.

\subsection{Proofs}

We first briefly prove for the one-dimensional case (i.e. Theorem
\ref{T4.1.1}), then extend the argument to high dimensions. It turns
out, comparing with the one-dimensional case, that the difficulty
point of the proof for high dimensions comes from the angle part.
So, a restriction concerning the angle part was made in Theorem
\ref{T4.1.2}.

\ \newline\emph{Proof of Theorem \ref{T4.1.1}.} Let $y(x):=
\Phi_\dd^{-1}\circ\varphi(x),\ x\ge\dd.$ We have

\beg{equation*}\beg{split} \ff{\d\mu}{\d y} &= \ff{\d\mu}{\d
x}\cdot \ff{\d x}{\d y}=
\e^{V(x)}\ff{\d\varphi^{-1}\circ\Phi_\dd(y)}{\d y}\\
& =\ff{\e^{V(x)}
\Phi_\dd'(y)}{\varphi'\circ\varphi^{-1}\circ\Phi_\dd(y)}
=\ff{\e^{V(x)}\Phi_\dd'(y)}{\varphi'(x)}=\Phi_\dd'(y).\end{split}\end{equation*}
Therefore, $\mu$ is the standard Gaussian measure under the new
coordinate $y\in [\dd,\infty)$. In other words, one has

$$\gamma(\d x):= (\mu\circ y^{-1})(\d x)= Z 1_{[\dd,\infty)}(x)\e^{-x^2/2}\d x,$$
where $Z$ is the normalization constant. By the HWI inequality
proved in \cite{BGL, OV, OV2} and the Gross log-Sobolev inequality
which implies the Talagrand inequality,  we have

\beq\label{4.2.1}\beg{split}&\gamma (g^2\log g^2) +W_2 (\gamma,
g^2\gamma)^2 \le 2
\ss{2\gamma((g')^2)}\, W_2(\gamma, g^2\gamma),\\
& W_2(\gamma, g^2\gamma)^2\le 2\gamma(g^2\log g^2),\ \ \
\gamma(g^2)=1.\end{split}\end{equation}  We remark
that although the HWI and Gross's log-Sobolev inequalities are
stated in the above references for the global Gaussian measure, they
are also true on a regular convex domain $\OO$, since the stronger
gradient estimate

$$|\nn P_t f|\le \e^{-t}P_t|\nn f|,\ \ \ f\in C_b^1(\OO)$$ holds for
the Neumann heat semigroup on $\OO$ (cf. \cite{W97} and references
within).

For any $f\in C_0^1([\dd,\infty))$ with $\mu(f^2)=1$, let $g:=
f\circ y^{-1}.$ We have

$$\ff{\d g}{\d x}= (f'\circ y^{-1})\ff{\d y^{-1}}{\d x}= \ff{f'\circ
y^{-1}}{y'\circ y^{-1}}= (f'\circ y^{-1})
\Big(\ff{\Phi_\dd'\circ\Phi_\dd^{-1}\circ\varphi}{\varphi'}\Big)\circ
y^{-1}.$$ Since $\gamma=\mu\circ y^{-1}$, this and (\ref{4.2.1})
imply (\ref{4.1.3}) and (\ref{4.1.4}). Finally, (\ref{4.1.2}) is implied
by (\ref{4.1.4}).\qed

\ \newline\emph{Proof of Theorem \ref{T4.1.2}.} Let $(r,\theta)$ be
the polar coordinate introduced in Section 2, and let $\nn_\theta$
denote the gradient operator on $\S^{d-1}$ for the standard metric
induced by the Euclidean metric on $\R^d$. By the orthogonal
decomposition of the gradient, we have

\beq\label{3.*} \nn f= (\pp_r f) \ff{\pp}{\pp r} + r^{-1}\nn_\theta
f,\ \ \ \ |\nn f|^2 = (\pp_r f)^2 + r^{-2} |\nn_\theta
f|^2.\end{equation} Let us introduce a new polar coordinate $(\bar
r, \theta)$, where

$$\bar r (r,\theta):= \Phi_0^{-1}\circ \varphi_\theta(r),\ \ \
r\ge 0, \theta\in \S^{d-1}.$$ We have

$$\d\mu := G(r,\theta)\d r\d\theta =\ff{G(r,\theta)}{\pp_r \bar
r}\d\bar r\d \theta= c(d) h(\theta) \Phi_0'(\bar r)\d\bar
r\d\theta=c(d) h(\theta)\d\mu_0,$$ where   $\d\mu_0:= \Phi_0'(\bar
r)\d\bar r\d\theta$ is the standard Gaussian measure under the new
polar coordinate $(\bar r,\theta)$. Thus, letting

$$y(x):= \bar r (x) \theta (x)= \Phi_0^{-1}\circ\varphi_{\ff{x}{|x|}}(|x|)\theta(x),
\ \ \ x\in\R^d,$$ we have

$$(\mu\circ y^{-1})(\d x)= c(d) h(x/|x|) (\mu_0\circ y^{-1})(\d x)
= c(d) h(x/|x|) \gamma(\d x),$$ where $\gamma$ is the standard
Gaussian measure on $\R^d$. By Gross' log-Sobolev inequality one has

$$\gamma(g^2\log g^2)\le 2 \gamma(|\nn g|^2),\ \ \ g\in C_0^\infty(\R^d), \mu_0(g^2)=1.$$
Thus, by the perturbation of the log-Sobolev inequality (cf.
\cite{DS2}), we have

\beq\label{4.2.2} (\mu\circ y^{-1})(g^2\log g^2)\le 2 C(h) (\mu\circ
y^{-1})(|\nn g|^2),\ \ \ g\in W^{2,1}(\gamma), (\mu\circ
y^{-1})(g^2)=1.\end{equation} Moreover, by \cite[Corollary
3.1]{BGL}, (\ref{4.2.2}) implies

\beq\label{4.2.3} W_2(\mu\circ y^{-1}, g^2\mu\circ y^{-1})^2\le 2
C(h) (\mu\circ y^{-1})(g^2\log g^2),\ \ \ (\mu\circ
y^{-1})(g^2)=1.\end{equation} This implies (\ref{4.1.3})
for the desired distance $\rr$ by using the change of variables
theorem as explained above.

Similarly, to prove  (\ref{4.1.2}) we intend apply (\ref{4.2.2}) for
$g:=f\circ y^{-1}$, where $f\in C_0^\infty(\R^d)$ with $\mu(f^2)=1.$
Since $y^{-1}= (\varphi_\theta^{-1}\circ\Phi_0(r), \theta)$ under
the polar coordinate, by the chain rule we have

$$\nn_\theta (f\circ y^{-1})=
\nn_\theta
f(\varphi_\theta^{-1}\circ\Phi_0(r),\theta)=\big((\nn_\theta f)\circ
y^{-1} +(\pp_r f)\circ y^{-1}\big)
\nn_\theta\varphi_\theta^{-1}\circ \Phi_0(r).$$ But
$\varphi_\theta\circ\varphi_\theta^{-1}\circ \Phi_0=\Phi_0$ implies

$$(\nn_\theta \varphi_\theta)(\varphi_\theta^{-1}\circ\Phi_0(r))
+ {\varphi_\theta}'\circ\varphi_\theta^{-1}\circ\Phi_0(r)\cdot
\nn_\theta(\varphi_\theta^{-1}\circ\Phi_0(r))=0,$$ where
$(\nn_\theta \varphi_\theta)(\varphi_\theta^{-1}\circ\Phi_0(r)):=
\nn_\theta
\varphi_\theta(s)|_{s=\varphi_\theta^{-1}\circ\Phi_0(r)},$ we arrive
at

\beg{equation}\label{3.*2}\beg{split} &|\nn_\theta (f\circ
y^{-1})|^2\\ &\le (1+\vv)(\pp_r f)^2\circ y^{-1}\Big(\ff{|\nn_\theta
\varphi_\theta(r)|(\varphi_\theta^{-1}\circ
\Phi_0(r))}{{\varphi_\theta}'\circ
\varphi_\theta^{-1}\circ\Phi_0(r)}\Big)^2+ (1+\vv^{-1})|\nn_\theta
f|^2\circ y^{-1}\\
&=(1+\vv)(\pp_r f)^2\circ y^{-1}\Big(\ff{|\nn_\theta
\varphi_\theta(r)|}{{\varphi_\theta}'(r)}\Big)^2\circ y^{-1}
+(1+\vv^{-1})|\nn_\theta f|^2\circ y^{-1}\end{split}\end{equation}
for any $\vv>0.$

On the other hand,
$$\pp_r (f\circ y^{-1}) = (\pp_r f)\circ y^{-1}
\ff{\Phi_0'(r)}{{\varphi_\theta}'\circ\varphi_\theta^{-1}\circ\Phi_0(r)}.$$
Since

\beq\label{3.*3} r=
\Phi_0^{-1}\circ\varphi_\theta(r(y^{-1}))=\Phi_0^{-1}\circ\varphi_\theta(r)\circ
y^{-1},\end{equation} we have

$$\Phi_0'(r)=
\big(\Phi_0'\circ\Phi_0^{-1}\circ\varphi_\theta(r)\big)\circ
y^{-1},\ \ \
{\varphi_\theta}'\circ\varphi_\theta^{-1}\circ\Phi_0(r)=
{\varphi_\theta}'(r)\circ y^{-1}.$$ Thus,

$$|\pp_r( f\circ y^{-1})|^2 = \Big\{(\pp_r
f)\ff{\Phi_0'\circ\Phi_0^{-1}\circ\varphi_\theta(r)}{{\varphi_\theta}'(r)}\Big\}^2\circ
y^{-1}.$$ Combining this with (\ref{3.*}), (\ref{3.*2}) and
(\ref{3.*3}), we obtain

\beg{equation*}\beg{split} &|\nn (f\circ y^{-1})|^2= (\pp_r
(f\circ
y^{-1}))^2 +r^{-2} |\nn_\theta (f\circ y^{-1})|^2\\
&\le \Big\{(\pp_r
f)\ff{\Phi_0'\circ\Phi_0^{-1}\circ\varphi_\theta(r)}{{\varphi_\theta}'(r)}\Big\}^2\circ
y^{-1}\\
&\quad +(\Phi_0^{-1}\circ\varphi_\theta(r))^{-2}\circ
y^{-1}\Big\{(1+\vv)(\pp_r f)^2\Big(\ff{|\nn_\theta \varphi_\theta(r)
|}{{\varphi_\theta}'(r)}\Big)^2+ (1+\vv^{-1})
|\nn_\theta f|^2\Big\}\circ y^{-1}\\
&= (\pp_r f)^2\circ
y^{-1}\Big\{\ff{(\Phi_0'\circ\Phi_0^{-1}\circ\varphi_\theta(r))^2}{({\varphi_\theta}'(r))^2}
+\ff{(1+\vv)|\nn_\theta \varphi_\theta(r)
|^2}{({\varphi_\theta}'(r))^2(\Phi_0^{-1}\circ
\varphi_\theta(r))^2}\Big\}\circ y^{-1}\\
&+\qquad (r\circ y^{-1})^{-2} |\nn_\theta f|^2\circ y^{-1}
\Big(\ff{(1+\vv^{-1})r^2}{(\Phi_0^{-1}\circ\varphi_\theta(r))^2}\Big)\circ
y^{-1}\\
 &\le |\nn f|^2\circ y^{-1}\max\Big\{\ff
{(1+\vv^{-1})r^2}{(\Phi_0^{-1}\circ\varphi_\theta
(r))^2},\ff{(\Phi_0'\circ\Phi_0^{-1}\circ\varphi_\theta(r))^2}{({\varphi_\theta}'(r))^2}
+\ff{(1+\vv)|\nn_\theta \varphi_\theta
(r)|^2}{({\varphi_\theta}'(r))^2(\Phi_0^{-1}\circ
\varphi_\theta(r))^2}\Big\}\circ y^{-1}\end{split}\end{equation*}
for any $\vv>0.$ Therefore,

\beq\label{4.2.4}|\nn (f\circ y^{-1})|^2\le(\aa |\nn f|^2)\circ
y^{-1}\end{equation} and hence (\ref{4.1.2})  follows from
(\ref{4.2.2}) by letting $g=f\circ y^{-1}.$

Finally, if $h$ is constant then $\mu\circ y^{-1}$ is the standard
Gaussian measure. Hence, by \cite[Theorem 4.3]{BGL} one has

\beg{equation*}\beg{split}&W_2(\mu\circ y^{-1}, (f^2\circ
y^{-1})\mu\circ y^{-1})^2 +(\mu\circ y^{-1})(f^2\circ y^{-1}\log
f^2\circ y^{-1})\\
&\le 2\ss{2(\mu\circ y^{-1})(|\nn (f\circ y^{-1})|^2)}\,
W_2(\mu\circ y^{-1}, (f^2\circ y^{-1}) \mu\circ
y^{-1}).\end{split}\end{equation*} By combining this with
(\ref{4.2.4}) we prove (\ref{4.1.7}). \qed

\ \newline\emph{Proof of Corollary \ref{C4.1.3}.} Since there exists
a constant $c_0>0$ such that

$$\Phi_0'(r)= c_0 r^{d-1}\e^{-r^2/2}=\beg{cases}
\Theta (r^{d-1}) &\text{as\ }r\to 0,\\
\Theta (r(1-\Phi_0(r))) &\text{as}\ r\to\infty,\end{cases}$$ where
$f=\Theta(g)$ means that the two positive functions $f$ and $g$ are
asymptotically bounded by each other up to constants, there exists a
constant $c\ge 1$ such that

$$\ff 1 c \Phi_0'(r)\le  \min\{r,r^{d-1}\} (1-\Phi_0(r))\le c\Phi_0'(r),\ \ \ r\ge 0.$$
Equivalently,

\beq\label{4.2.5} \ff 1 c \Phi_0'\circ\Phi_0^{-1}(r)\le
\min\{\Phi_0^{-1}(r),\ \Phi_0^{-1}(r)^{d-1}\}(1-r)\le
c\Phi_0'\circ\Phi_0^{-1}(r),\ \ \ r\in [0,1).\end{equation}

Next, it is easy to see from (\ref{V}) that

\beq\label{4.2.6}\Phi_0^{-1}\circ\varphi_\theta (r) = \beg{cases}
\Theta(r^{\dd/2})\ &\text{as}\ r\to\infty,\\
\Theta(r)\ &\text{as}\ r\to 0,\end{cases}
\end{equation} and

\beq\label{4.2.7}\ff{1-\varphi_\theta(r)}{{\varphi_\theta}'(r)}
=\ff{\int_r^\infty s^{d-1}\e^{V(s\theta)}\d s}{r^{d-1}
\e^{V(r\theta)}}\le c r^{1-\dd}\end{equation} for some constant
$c>0$ and all $r\ge 1.$ Combining (\ref{4.2.5}), (\ref{4.2.6}) and
(\ref{4.2.7}) we obtain

\beq\label{MM}\max\Big\{\ff{r^2}{(\Phi_0^{-1}\circ
\varphi_\theta(r))^2},\ \ff{(\Phi_0'\circ
\Phi_0^{-1}\circ\varphi_\theta(r))^2}{({\varphi_\theta}'(r))^2}\Big\}\le
c (1+r)^{2-\dd}\end{equation} for some constant $c>0$.

If (\ref{V2}) holds then

$$|\nn_\theta \varphi_\theta(r)|=|\nn_\theta(1-\varphi_\theta(r)| \le c_4 \min\Big\{r^d, \int_r^\infty s^{d-1}
\e^{V(s\theta)}\d s\Big\},$$ so that due to (\ref{4.2.6}) and
(\ref{4.2.7})

$$\ff{|\nn_\theta\varphi_\theta(r)|^2}{({\varphi_\theta}'(r))^2(\Phi_0^{-1}\circ\varphi_\theta(r))^2} \le
c_5\bigg(\ff{\min \{r^d, \int_r^\infty s^{d-1} \e^{V(s\theta)}\d s
\}}{(r1_{\{r<1\}}+ r^{\dd/2}1_{\{r\ge 1\}})r^{d-1}
\e^{V(r\theta)}}\bigg)^2\le c_6 (1+r)^{2-3\dd}$$ for some constants
$c_5,c_6>0$. Combining this with (\ref{MM}) and Theorem \ref{T4.1.2},
we prove (\ref{4.1.9}).

Finally, for any $x_1,x_2\in \R^d$ let $i\in \{1,2\}$ such that
$|x_i|=|x_1|\lor |x_2|.$ Similarly to the proof of Corollary
\ref{C1.2}, define

$$f(x)= \ff{|x-x_i|\land \ff{|x_i|} 2 }{(1+|x_i|)^{1-\dd/2}},\ \ \ x\in \R^d.$$ Then

$$\GG(f,f):= (1+|\cdot|)^{2-\dd}|\nn f|^2\le \ff {1_{\{|x_i|/2\le |\cdot|\le 3|x_i|/2\}}(1+|\cdot|)^{2-\dd}} {(1+|x_i|)^{2-\dd}}
\le C(\dd)$$ for some constant $C(\dd)>0.$  Since $|x_i|\ge \ff 1 2
|x_1-x_2|$, this implies that the intrinsic distance $\rr$  induced
by $\GG$ satisfies

$$\rr(x_1,x_2)^2 \ge \ff {|f(x_1)-f(x_2)|^2} {C(\dd)} \ge
C_1(\dd)\tt\rr(x_1,x_2)^2$$ for some constant $C_1(\dd)>0$, and
hence is complete. Thus, by \cite[Theorem 1.1]{W04} or \cite[Theorem
6.3.3]{Wbook}, (\ref{4.1.10}) follows from (\ref{4.1.9}).\qed

\paragraph{Acknowledgement.} The author would like to thank the
referees for useful comments.

\beg{thebibliography}{99}


\bibitem{BLW} D. Bakry, M. Ledoux and F.-Y. Wang,
\emph{Perturbations of functional inequalities using growth
conditions,}   J. Math. Pures Appl. 87(2007), 394-407.


\bibitem{BGL} S. G. Bobkov, I.  Gentil and M.  Ledoux,
\emph{Hypercontractivity of Hamilton-Jacobi equations,} J. Math.
Pures Appl. 80(2001), 669--696.

\bibitem{BG} S. G. Bobkov, and F. G\"otze, \emph{Exponential integrability and
transportation cost related to logarithmic Sobolev inequalities,} J.
Funct. Anal. 163(1999), 1--28.

\bibitem{BV} F. Bolley and C. Villani, \emph{Weighted
csisz\'ar-Kullback-Pinsker inequalities and applications to
transportation inequalities,} Ann. Fac. Sci. Toulouse Math. (6),
14(2005), 331--352.

\bibitem{Ca} L. A. Caffarelli, \emph{The regularity of mappings
with a convex potential,} J. Amer. Math. Soc. 5(1992), 99--104.


\bibitem{DS} J. D. Deuschel and D. W. Stroock, \emph{Large Deviations,}
Pure and Appl. Math. Ser. 137, Academic Press, San Diego, 1989.

\bibitem{DS2} J. D. Deuschel and D. W. Stroock,
\emph{Hypercontractivity and spectral gap of symmetric diffusions
with applications to the stochastic Ising models,} J. Funct. Anal.
92(1990), 30--48.

\bibitem{DGW} H. Djellout, A. Guilin and L.-M. Wu,
\emph{Transportation cost-information inequalities for random
dynamical systems and diffusions,} Ann. Probab. 32(2004),
2702--2732.

\bibitem{GGM} I. Gentil, A. Guillin and L. Miclo, \emph{Modified
logarithmic Sobolev inequalities and transportation inequalities,}
Probab. Theory Relat. Fields 133(2005), 409--436.

\bibitem{GGM2} I. Gentil, A. Guilin and L. Miclo, \emph{Modified
logarithmic Sobolev inequalities in null curvature,}  Revista Mat.
Ibero. 23(2007).

\bibitem{Go} N. Gozlan, \emph{Characterization of Talagrand's like
transportation-cost inequalities on the real line,}   J.
Funct. Anal. 250(2007), 400--425.

\bibitem{Gross} L. Gross,  \emph{Logarithmic Sobolev inequalities,} Amer. J. Math.
{  97}(1976), 1061--1083.

\bibitem{GW} F.-Z.  Gong and F.-Y. Wang, \emph{
 Functional inequalities for uniformly integrable semigroups and
application to essential spectrums,} Forum Math.  14(2002),
293--313.
\bibitem{G06} N. Gozlan, \emph{Integral criteria for transportation
cost inequalities,} Electron. Comm. Probab. 11(2006), 64--77.

\bibitem{K} A. V. Kolesnikov, \emph{Convexity inequalities and
optimal transport of infinite dimensional measures,} J. Math. Pures
Appl. 83(2004), 1373--1404.

\bibitem{Ma} K. Marton, \emph{A simple proof of the blowing-up
lemma,} IEEE Trans. Inform. Theory, 32(1986), 445--446.

\bibitem{OV} F.  Otto and C. Villani, \emph{Generalization of an inequality
by Talagrand and links with the logarithmic Sobolev inequality,}
J. Funct. Anal.  173(2000), 361--400.

\bibitem{OV2} F.  Otto and C. Villani, \emph{Comment on: ``Hypercontractivity of Hamilton-Jacobi equations",
by S. Bobkov, I. Gentil and M. Ledoux,} J. Math. Pures Appl.
80(2001), 697--700.

\bibitem{RW01} M. R\"ockner and F.-Y. Wang,   \emph{Weak Poincar\'e inequalities
and $L^2$-convergence rates of Markov semigroups,} J. Funct. Anal.
185(2001), 564--603.


\bibitem{Shao} J. Shao, \emph{Hamilton-Jacobi semigroup in infinite
dimensional spaces}, Bull. Sci. Math. 130(2006), 720--738.

\bibitem{T} M. Talagrand, \emph{Transportation cost for Gaussian
and other product measures,} Geom. Funct. Anal. 6(1996), 587--600.

\bibitem{W97}    F.-Y. Wang, \emph{Logarithmic {S}obolev inequalities on noncompact
  {R}iemannian manifolds}, Probab. Theory Related Fields  109 (1997), 417--424.

\bibitem{W00} F.-Y. Wang, \emph{Functional inequalities for empty essential
spectrum,} J. Funct. Anal. 170(2000), 219--245.

\bibitem{W00b} F.-Y. Wang, \emph{Functional inequalities, semigroup properties
and spectrum estimates,} Infin. Dimens. Anal. Quant. Probab.
Relat. Topics 3(2000), 263--295.


\bibitem{W04} F.-Y. Wang,  \emph{Probability distance
inequalities on Riemannian manifolds and path spaces,} J. Funct.
Anal. 206(2004), 167--190.


\bibitem{Wbook} F.-Y. Wang, \emph{Functional Inequalities, Markov Processes, and Spectral Theory}, Science
Press, Beijing 2005.


\end{thebibliography}

\end{document}